\documentclass{amsart}
\usepackage{amssymb, amsmath, latexsym, amsthm, nccmath}
\usepackage{fancyhdr}
\usepackage{graphicx}
\usepackage[margin=1in]{geometry}
\usepackage{multicol}

\newtheoremstyle{dotless}{}{}{\itshape}{}{\bfseries}{}{ }{}
\theoremstyle{dotless}

\newtheorem{theorem}{Theorem}[section]
\newtheorem{lemma}{Lemma}%[section]
\theoremstyle{definition}

\newcommand{\La}{\langle}
\newcommand{\Ra}{\rangle}

\newcommand{\J}{J_0\vee s(J_0)}
\newcommand{\cD}{\mathcal D}

\newcommand{\bR}{\mathbb R}

\newcommand{\bfT}{\mathbf T}

\newcommand{\one}{\mathbf{1}}
\begin{document}

\title[Dyadic bi-parameter simple commutator and dyadic little BMO]%
{Dyadic bi-parameter simple  commutator and dyadic little BMO}
%\author[N.~Arcozzi]{Nicola Arcozzi}
%\address{Universit\`{a} di Bologna, Department of Mathematics, Piazza di Porta S. Donato, 40126 Bologna (BO)}
%\email{nicola.arcozzi@unibo.it}
%\thanks{Theorem 3.1 was obtained in the frameworks of the project 17-11-01064 by the Russian Science Foundation}
%\thanks{NA is partially supported by the grants INDAM-GNAMPA 10017 "Operatori e disuguaglianze integrali in spazi con simmetrie" and PRIN 10018 "Variet\`{a} reali e complesse: geometria, topologia e analisi armonica"}
\author[I.~ Holmes]{Irina Holmes}
\address{Department of Mathematics, Texas A\&M University}
\thanks{IH is partially supported by the NSF an NSF Postdoc under Award No.1606270}
%\author[P. Mozolyako]{Pavel Mozolyako}
%\thanks{PM is supported by the Russian Science Foundation grant 17-11-01064}
%\address{Universit\`{a} di Bologna, Department of Mathematics, Piazza di Porta S. Donato, 40126 Bologna (BO)}
%\email{pavel.mozolyako@unibo.it}
\author[S.~Treil]{Sergei Treil}
\thanks{ST is partially supported by the NSF  DMS 1856719}
\address{Department of Mathematics, Brown University, 151 Thayer Str./Box 1917, Providence, RI 02912, USA}
\email{treil@math.brown.edu \textrm{(S. Treil)}}
\author[A.~Volberg]{Alexander Volberg}
\thanks{AV is partially supported by the NSF grant DMS 1900268 and by Alexander von Humboldt foundation}
\address{Department of Mathematics, Michigan Sate University, East Lansing, MI. 48823}
\email{volberg@math.msu.edu \textrm{(A.\ Volberg)}}
\makeatletter
\@namedef{subjclassname@10010}{
  \textup{10010} Mathematics Subject Classification}
\makeatother
\subjclass[10010]{42B100, 42B35, 47A100}
% 42B	Harmonic analysis in several variables
% 42B100	Singular and oscillatory integrals (Calder?on-Zygmund, etc.)
% 42B35	Function spaces arising in harmonic analysis
% 47A	General theory of linear operators
% 47A100	Norms (inequalities, more than one norm, etc.)
%{100E100, 47B37, 47B40, 100D55.}
%
% 100D55	$H^p$-classes (1980-10009)
% 100E100	Integration, integrals of Cauchy type, integral representations of analytic functions
%
% 47B   	Special classes of linear operators
% 47B37	Operators on special spaces (weighted shifts, operators on sequence spaces, etc.)
% 47B40	Spectral operators, decomposable operators, well-bounded operators, etc.
\keywords{Bi-parameter Carleson embedding, bi-parameter Hankel operators, dyadic model}
\begin{abstract} 
Let $\bfT$ is a certain tensor product of simple dyadic shifts defined below.
We prove here that for dyadic bi-parameter commutator the following equivalence holds
$  \|\bfT  b-b \bfT \| \asymp \|b\|_{bmo^d}$.  This result is well-known for many types of bi-parameter commutators, see \cite{FS}, \cite{DLWY} and \cite{DPSK} for more details.
\end{abstract}
\maketitle

\section{Introduction}
\label{intro}

In this note we are considering a simple dyadic model of bi-parameter commutator. Bi-parameter theory is notoriously more difficult than the more classical one  parameter theory of singular integrals. The good place to get acquinted with multi-parameter specifics are the papers of J.-L. Journ\'e \cite{JLJ}, \cite{JLJ2} and Muscalu--Pipher--Tao--Thiele \cite{MPTT1}, \cite{MPTT2}. The applications to analysis in polydisc can be found in \cite{Ch},  \cite{ChF}, \cite{ChF2}.

It is well known \cite{ChF, ChF2, JLJ, JLJ2} that in the multi-parameter setting all concepts of Carleson measure, $BMO$, John--Nirenberg inequality, Calder\'on--Zygmund decomposition (used in classical theory) are much more delicate. Paper \cite{MPTT1}  develops a completely new approach to prove natural tri-linear bi-parameter estimates on bi-parameter paraproducts, especially outside of Banach range. In \cite{MPTT1} Journ\'e's lemma \cite{JLJ2} was used, but the approach did not generalize to multi-parameter paraproduct forms. This issue was resolved in \cite{MPTT2}, where a simplified method was used  to address the multi-parameter paraproducts.

One of the new feature of the multi-parameter theory is captured by several different definitions of $BMO$, see \cite{FS}, \cite{Ch}. The necessity of those new effects was discovered first by Carleson \cite{Car}, see also \cite{Tao}. The difficulties with multi-parameter theory was highlighted recently by two very different series of papers. One concerns with the Carleson measure and Carleson embedding on polydisc and on multi-tree, see, e.g. \cite{AMPS18}, \cite{AHMV}, \cite{MPV}, \cite{MPVZ}.  These papers, roughly speaking, are devoted to harmonic analysis (Carleson embedding in particular) on graphs with cycles. This theory is drastically different from the usual one, and it is much more difficult. Another particularity of the multi-parameter theory is highlighted by the repeated commutator story. Lacey and Ferguson \cite {FL} found the characterization of the symbols that give us bounded ``small Hankel operators". It was a breakthrough article that gave a multi-parameter Nehari theorem and  a long searched after factorization of bi-parameter Hardy space $H^1$. It was exactly equivalent to a bi-parameter ``repeated commutator characterization". Several papers followed where ``bi-'' was upgraded to ``multi-'', and where repeated commutation was performed with different classical singular integrals: \cite{DP}, \cite {FL}, \cite{L}, \cite {LPPW}, \cite {LT}.  

There are many new and beautiful ideas in the above mentioned papers devoted to this subject. We also want to mention  \cite{OPS}, where the authors use an argument inspired by Toeplitz operators to  
show the lower bounds. It assumes the lower norm for the Hilbert  
transform  claimed in \cite{FL}, \cite{LPPW} as a black  
box.  But there is a problem with \cite{FL} and with what followed.

That was a big breakthrough in multi-parameter theory. Unfortunately \cite{V} indicated a hole in all the proofs of \cite{FL}, \cite{L}, and this circle of problems is still unsettled.

We decided to simplify the problem to the ``bare bones", and to consider commutators with the simples dyadic singular integral--the dyadic shift of order $1$. The repeated commutator problem presented the same type of difficulty as in \cite{FL} and subsequent papers. In our mind, it is very nice to see the difficulty in a dyadic case removed from all technicalities of a continuous case. But here  we consider a much simpler problem: the characterization of the boundedness of the usual (not repeated) commutator of multiplication and the tensor product of two dyadic singular integrals. Again, the dyadic singular integrals are the simplest possible, they are the simplest dyadic shifts.

Not surprisingly we are able to solve this problem completely and to get an expected characterization of the symbol of the bounded simple commutator: it is small (dyadic) $bmo$. This answer is not surprising at all, as the similar results (for various singular operators) were obtained in \cite{FS},  \cite{DLWY}, \cite{P}, \cite{DPSK}.

Notice that for repeated commutator {\it the expected  characterization} of the symbol of the bounded repeated commutator is a different $BMO$, namely, it
is $BMO_{ChF}$: this is Chang--Fefferman $BMO$ (or product $BMO$) studied in \cite{ChF}, \cite{ChF2}. The fact that 
$BMO_{ChF}$ is an expected characterization of bounded repeated commutators has at least two ``confirmations". One is the paper of Blasco--Pott \cite{BP}, where it is proved that the dyadic $BMO_{ChF}$ characterizes  the boundedness of repeated commutator with {\it all dyadic martingale transforms simultaneously}. What one would like to show, that just one interesting transform (dyadic) is enough. The second ``confirmation" is that the counterexample found in \cite{V} seems to be ``easily circumvented''. Unfortunately the attempt to do that did not suceed--to the best of our knowledge--so far.

Let us finish by writing a simple commutator (we deal with it below) and repeated commutator (that brings so much pain). Let $\bfT= T_1\otimes T_2$, where $T_1, T_2$ are two (dyadic) singular integrals, each acting in its own $L^2(\bR)$: $T_1$ acts on functions of variable $x$, $T_2$ acts on functions whose variable is called $y$. Now let $b(x, y)$ be a symbol. We need to characterize the boundedness of simple commutator $[\bfT, b] =T_2T_1 b- bT_2T_1$ in terms of $b$.  We do this below. On the other hand, if one is concerned with repeated commutator, then one is considering the following nested commutation: $[T_2, [T_1, b]]=T_2(T_1 b- bT_1)-(T_1 b- bT_1)T_2= T_2T_1 b+ bT_1T_2 - T_2 b T_1 - T_1 bT_2 $.

\section{Plan}
\label{plan}

Let $\cD$ be the usual dyadic lattice on the line and $\cD\times \cD$ be the family of dyadic rectangles on the plane.
Haar functions are denoted $\{h_I\}_{I \in \cD}$. We consider $R=I\times J$ and $h_R= h_I\otimes h_J$.
Let $T$ be the dyadic shift defined by
$$
T h_I = h_{I_+}- h_{I_-}, \,\, \text{if}\,\, I\,\, \text{is even}; \quad Th_I=0,\,\, \text{otherwise},
$$
where we say $I \in \cD$ is even if $|I|=2^{-2k}$.

Let $\bfT := T\otimes T$.
We are interested in the commutator
$$
\bfT  b -b\bfT 
$$
and in characterization of its boundedness. The ultimate goal is to prove that its boundedness is equivalent to
$b\in bmo^d$, 
where $bmo^d$ is dyadic little $bmo$.

%-----------------------------------------------%
%---------Section: One-Param Commut-------------%
%-----------------------------------------------%
\section{The $1$-parameter case}
\label{D-part}

\begin{theorem}
\label{T:1d}
$[T, b]$ is $L^2$-bounded if and only if $b \in BMO^d$.
\end{theorem}

Recall that 
	$$\|b\|_{BMO^d}^2 = \sup_{I\in\cD} \frac{1}{|I|}\int_I |b(x) - \La b\Ra_I|\,dx =
		\sup_{I\in\cD} \frac{1}{|I|}\sum_{I' \subseteq I} |(b, h_I)|^2.$$
We will make use of the paraproduct decomposition
	$$b(x)f(x) = \pi(b, f)(x) + Z(b, f)(x) + D(b, f)(x),$$
where
	\begin{equation}
	\label{E:para-def}
	\pi(b, f) := \sum_{I} (b, h_I) \La f\Ra_I h_I,
	\:\:\:\:\:
	Z(b, f) := \sum_{I} (b, h_I)(f, h_I)\tilde\one_I,
	\:\:\:\:\:
	D(b, f) := \sum_{I} \La b\Ra_I (f, h_I) h_I.
	\end{equation}
Above, summation is always over $I \in \cD$, and we used the recurring notations
	$$ (f, h_I) := \int f h_I\,dx; 
	\:\:\:\:\:\:
	\La f\Ra_I := \frac{1}{|I|}\int_I f(x)\,dx;
	\:\:\:\:\:\:
	\tilde\one_I := \frac{\one_I}{|I|}.$$

Key to the proof will be the following:
\begin{lemma}
\label{dia}
If $[T, b]$ is $L^2$-bounded, then 
	\begin{equation}
	\label{onestep}
	|(b, h_I)|^2 \le C_0|I|
	\end{equation} for all intervals $I\in \cD$.
\end{lemma}
We prove this result below in section \ref{sS:proofdia}.

%------lower bound 1par------%
\subsection{Proof of Theorem \ref{T:1d} -- Lower bound.}
\begin{proof}[\unskip\nopunct]

This part shows that if $[T, b]$ is bounded then $b \in BMO^d$.
We split the $Tb-bT$ operator into $D, \pi, Z$ pieces:
	$$
	[T, b]f = T(bf) - b \: T f =
	\left[T\pi(b, f) - \pi(b, Tf)\right] +
	\left[TD(b, f) - D(b, Tf)\right] +
	\left[TZ(b, f) - Z(b, Tf)\right].
	$$
We consider the function $b$ to be fixed throughout, so we refer to these operators simply as $(T\pi-\pi T)$, 
$TD - DT$, and $TZ - ZT$.

\vspace{0.1in}
\noindent \textbf{1. $D$ part.} We have
	$$
	(TD - DT)f = \sum_{I even} \left( D(b, f), h_I\right) (h_{I_+} - h_{I-}) - \sum_I \La b\Ra_I (Tf, h_I) h_I.
	$$
Now, by definition of $T$,
	$$(Tf, h_I) = \left\{ \begin{array}{rl}
	(f, h_{\hat I}) s(I, \hat I), & \text{ if } I \text{ odd}\\
	0, & \text{ otherwise},
	\end{array}\right.$$
where, for every $I \in \cD$, $\hat I$ denotes the dyadic parent of $I$ and
	$$ s(I, \hat I) := \left\{ \begin{array}{rl}
		1, & \text{ if } I = \hat I_{+}\\
		-1, & \text{ if } I = \hat I_-.
	\end{array}\right.
	$$
So
	$$
	(TD - DT)f = \sum_{I even} \La b\Ra_I (f, h_I) (h_{I_+} - h_{I-}) - 
			\sum_{I odd} \La b\Ra_I (f, h_{\hat I}) s(I, \hat I) h_I
	$$
Relabeling the first term over the odd intervals instead,
	\begin{eqnarray*}
	(TD - DT)f &=& \sum_{I odd} \La b\Ra_{\hat I} (f, h_{\hat I}) s(I, \hat I) h_I - 
			\sum_{I odd} \La b\Ra_I (f, h_{\hat I}) s(I, \hat I) h_I \\
			&=& \sum_{I odd} \left(\La b\Ra_{\hat I} - \La b\Ra_I\right) (f, h_{\hat I}) s(I, \hat I) h_I.
	\end{eqnarray*}
It is easy to see that
	$$
	\La b\Ra_{\hat I} - \La b\Ra_I = \frac{-s(I, \hat I)}{\sqrt{|\hat I|}} (b, h_{\hat I}),
	$$	
so
	$$
	(TD - DT)f = \sum_{I odd} \frac{-1}{\sqrt{|\hat I|}} (b, h_{\hat I}) (f, h_{\hat I}) h_I.
	$$
This convenient formula for this term, combined with Lemma \ref{dia}, shows that $TD-DT$ is actually \textit{a priori bounded} on all test functions: 
	$$
	\|(TD - DT)f\|_2^2 = \sum_{I odd} \frac{1}{|\hat I|} |(b, h_{\hat I})|^2 |(f, h_{\hat I})|^2 
		\leq C_0 \sum_{I odd} |(f, h_{\hat I})|^2  \leq \|f\|_2^2.
	$$
Therefore we can forget about this term completely, and focus on the remaining two.

\vspace{0.1in}
\noindent \textbf{2. $Z$ part.} This term does not simplify as nicely as the $D$-term:
	$$
	(TZ - ZT)f = \sum_I (b, h_I) (f, h_I) T(\tilde \one_I) -
		\sum_I (b, h_I) (Tf, h_I) \tilde\one_I.
	$$
So now we test our commutator on a Haar function $h_{I_0}$, where $I_0$ is {\it even}. In particular, the $Z$-term gives
	$$
	(TZ - ZT)h_{I_0} = (b, h_{I_0}) T\tilde\one_{I_0} - (b, h_{I_0^+}) \tilde\one_{I_0^+}
		+ (b, h_{I_0^-}) \tilde\one_{I_0^-}.
	$$
From Lemma \ref{dia}, we see immediately that the $(TZ-ZT)$ term is a priori uniformly bounded on $h_{I_0}$, so it can be put aside.

\vspace{0.1in}
\noindent \textbf{3. $\pi$ part.}
Now we are left with
	$$
	(T\pi - \pi T)f = \sum_{I even} (b, h_I) \La f\Ra_I (Th_I) - \sum_I (b, h_I) \La Tf\Ra_I h_I,
	$$
which we write as
	$$(T\pi - \pi T)f = \Sigma_1 + \Sigma_2,$$
where
	\begin{eqnarray*}
	\Sigma_1 &=& - \sum_{I even} (b, h_I) \La Tf\Ra_I h_I\\
	\Sigma_2 &=& \sum_{I even} (b, h_I) \La f\Ra_I (Th_I) - \sum_{I odd} (b, h_I) \La Tf\Ra_I h_I.
	\end{eqnarray*}	
Now, note that $\Sigma_1$ and $\Sigma_2$ are \textit{orthogonal}, since $\Sigma_1$ lives in the {\it even} space and $\Sigma_2$ lives in the {\it odd} space. So, $(T\pi - \pi T)$ being uniformly bounded on $h_0$ means that the terms $\Sigma_1$ and $\Sigma_2$ are individually uniformly bounded on $h_0$.

So let us first plug in our test function into $\Sigma_1$:
$$
\Sigma_1 h_{I_0} = \sum_{\substack{I even \\ I \subsetneq I_0^+}} (b, h_I) h_{I_0^+}(I) h_I
	- \sum_{\substack{I even \\ I \subsetneq I_0^-}} (b, h_I) h_{I_0^-}(I) h_I,
	$$
where for dyadic intervals $I \subsetneq J$, $h_J(I)$ denotes the constant value that $h_J$ takes on the interval $I$. Then
	$$
	\|\Sigma_1 h_{I_0}\|_2^2 = \frac{2}{|I_0|} \sum_{\substack{I even\\ I \subsetneq I_0^\pm}} |(b, h_I)|^2 \lesssim C
	$$
is uniformly bounded. Combined with \eqref{onestep}, we get
	\begin{equation}\label{E:evenints1}
	\sum_{I \subset I_0, I even} |(b, h_I)|^2 \lesssim C|I_0|,
	\end{equation}
for all dyadic intervals $I_0$. Technically, so far this is only true for even $I_0$, but
if $I_0$ is odd, the same quickly follows by splitting into the subintervals of the even halves $I_0^\pm$ and then using \eqref{onestep} once more.

Since \eqref{E:evenints1} is only half the battle, we must look at $\Sigma_2$ now, which is also uniformly bounded on our test functions. We will split $\Sigma_2$ into its own two pieces:
	$$\Sigma_2 = \Sigma_{21} - \Sigma_{22},$$
where
	$$
	\Sigma_{21}f = \sum_{I even} (b, h_I) \La f\Ra_I (Th_I) \text{ and } 
	\Sigma_{22} = \sum_{I odd} (b, h_I) \La Tf\Ra_I h_I.
	$$
For $f = h_{I_0}$,
	$$
	\Sigma_{21}h_{I_0} = \sum_{\substack{I odd \\ \hat I \subsetneq I_0^\pm}} (b, h_{\hat I}) h_{I_0}(\hat I) s(I, \hat I) h_I,
	$$
so
	$$
	\|\Sigma_{21}h_{I_0}\|_2^2 = \frac{1}{|I_0|} \sum_{\substack{I odd \\ \hat I \subsetneq I_0^\pm}} |(b, h_{\hat I})|^2
		= \frac{1}{|I_0|} \sum_{\substack{I even \\   I \subsetneq I_0^\pm}} |(b, h_I)|^2,
		$$
which we know to be uniformly bounded already from \eqref{E:evenints1}. This means that the last remaining term, $\Sigma_{22}$ is again uniformly bounded on our test functions:
	$$\Sigma_{22}h_{I_0} = \sum_{\substack{I odd \\ I \subsetneq I_0^+}} (b, h_I) h_{I_0^+}(I) h_I 
		- \sum_{\substack{I odd \\ I \subsetneq I_0^-}} (b, h_I) h_{I_0^-}(I) h_I. $$
So
	$$
	\|\Sigma_{22}h_{I_0}\|_2^2 = \frac{2}{|I_0|} \sum_{\substack{I odd\\ I \subsetneq I_0^\pm}} |(b, h_I)|^2
	$$
is also uniformly bounded. Then the odd counterpart of \eqref{E:evenints1}:
	$$
	\sum_{I \subset I_0, I odd} |(b, h_I)|^2 = \sum_{\substack{I odd\\ I \subsetneq I_0^\pm}} |(b, h_I)|^2 + |(b, h_{I_0^+})|^2
	+ |(b, h_{I_0^-})|^2 \lesssim C|I_0|
	$$
follows again by combination with \eqref{onestep}.

So now we have 
	$$\sum_{I\subset I_0} |(b, h_I)|^2 \lesssim C |I_0|$$
for all $I_0 \in \cD$, or that $b \in BMO^d$. 
\end{proof}

%------base lemma proof 1par------%

\subsection{Proof of Lemma \ref{dia}}
\label{sS:proofdia}
\begin{proof}[\unskip\nopunct]
The boundedness of $Tb-bT$ implies uniform boundedness of 
	$$((Tb-bT)h_I,  h_K ) \text{ for all } I, K \in \mathcal{D}.$$
Now
	\begin{equation}\label{E:onestep1}
	\left( (Tb-bT)h_I, h_K\right) = (bh_I, \underbrace{T^* h_K}_{0 \text{ if } K \text{ even}}) - 
		(b \underbrace{Th_I}_{0 \text{ if } I \text{ odd}}, h_K),
	\end{equation}
since
	$$T^*h_J = \left\{ \begin{array}{ll}
		s(J, \hat J) h_{\hat J}, & \text{ if } J \text{ is odd};\\
		0, & \text{ if } J \text{ is even.}
	\end{array}\right.$$

First, in \eqref{E:onestep1}, take $I$ to be even and $K = I_\pm$ (a dyadic child of $I$). 
Then $T^*h_K = s(K, I) h_I$, so
	$$
	((Tb-bT)h_I, h_{I_{+}}) = (bh_I, s(K, I) h_I) - (b(h_{I_{+}} - h_{I_-}), h_K) 
		= s(K, I) \big(\La b\Ra_I - \La b\Ra_{K}\big) = \frac{\pm 1}{\sqrt{|I|}} (b, h_I), 
	$$
which proves \eqref{onestep} for all even $I$. Now take in \eqref{E:onestep1} an odd $I$ and $K = I^{(2)}$, the dyadic \textit{grandparent} of $I$ which will then also be odd. In this case, we have
	$$
	((Tb-bT)h_I, h_{I^{(2)}}) = (bh_I, T^* h_{I^{(2)}}) = (bh_I, \pm h_{I^{(3)}}) = \pm \frac{1}{\sqrt{8|I|}} (b, h_I),
	$$
and this proves \eqref{onestep} for all odd $I$.
\end{proof}

%------upper bound 1par------%
\subsection{Proof of Theorem \ref{T:1d} -- Upper bound.}
\begin{proof}[\unskip\nopunct]
The proof that $[T, b]$ is bounded if $b \in BMO^d$ is completely standard and well-known, but we include it here for completeness.
We show that each of the terms $(TD-DT)$, $(T\pi-\pi T)$ and $(TZ - ZT)$ is bounded. The easiest term is $(TD-DT)$, whose $L^2$-norm can just be computed to be
	$$
	\|(TD-DT)\|_2^2 = \sum_{I odd} \frac{1}{|\hat I|} |(b, h_{\hat I})|^2 |(f, h_{\hat I})|^2 \lesssim \|b\|_{BMO^d} \|f\|_2^2,
	$$
since obviously $|(b, h_I)| \leq \|b\|_{BMO^d} |I|$ for all $I$. For the other terms we can use $H_d^1 - BMO^d$ duality, which gives us that
	$$
	|(b, \Phi)| \lesssim \|b\|_{BMO^d}\|S_d\Phi\|_1,
	$$
where $S_d$ is the dyadic square function:
	$$
	S_d^2 f := \sum_{I} |(f, h_I)|^2\tilde\one_I.
	$$
	
For the $(T\pi-\pi T)$ term:
	$(T\pi-\pi T) f = \Pi_1 f - \Pi_2 f$, where
	$$
	\Pi_1 f = \sum_{I even} (b, h_I) \La f\Ra_I (Th_I); \:\:\:\:
	\Pi_2 f = \sum_I (b, h_I) \La Tf\Ra_I h_I.
	$$
Then for $f, g \in L^2(\bR)$, we have $|(\Pi_1f, g)| = |(b, \Phi_1| \lesssim \|b\|_{BMO^d} \|S_d\Phi_1\|_1$, where
	$$
	\Phi_1 = \sum_{I even} \La f\Ra_I (T^*g, h_I) h_I,
	$$
so
	$$
	S_d^2\Phi_1 = \sum_{I even} |\La f\Ra_I|^2 |(T^*g, h_I)|^2 \tilde\one_I \leq (M_d^2 f) (S_d^2T^*g),
	$$
where $M_d$ is the dyadic maximal function. Then (since $M^d$, $S^d$ and $T$ are bounded), $\|S_d\Phi_1\|_1 \lesssim \|f\|_2\|g\|_2$, and
$|(\Pi_1f, g)| \lesssim \|b\|_{BMO^d}\|f\|_2\|g\|_2$ for all $f, g \in L^2(\bR)$. So we have bounded the first term $\Pi_1$. Similarly,
$|(\Pi_2f, g)| = |(b, \Phi_2)|$ where
	$$
	\Phi_2 = \sum_I \La Tf\Ra_I (g, h_I) h_I.
	$$
Then
	$$
	S_d^2\Phi_2 = \sum_I \La Tf\Ra_I^2 |(g, h_I)|^2\tilde\one_I \leq M_d^2(Tf) S_d^2g,
	$$
giving again $|(\Pi_2f, g)| \lesssim \|b\|_{BMO^d} \|S_d\Phi_2\|_1 \lesssim \|b\|_{BMO^d} \|f\|_2 \|g\|_2$.

The term $(TD-DT)$ follows very similarly.
\end{proof}

%-----------------------------------------------%
%----------Section: Bi-Param Commut-------------%
%-----------------------------------------------%

\section{The bi-parameter commutator}
\label{bi-param-rec}

Now let us work in $\bR^2 = \bR\otimes\bR$ and consider $[\bfT, b]$, where $\bfT = T \otimes T$ and $b(x, y)$ is an $\bR^2$-function. We will denote $\bfT = T_1 \otimes T_2$, where $T_i$ means $T$ acting only on the $i^{\text{th}}$ variable. Dyadic intervals in first coordinate will be called only by letters $I, I', I'', K$, in the second coordinate by $J, J', J'', L$. Dyadic rectangles $R\in \cD^2$ are of the form $R = I\times J$, with $I, J \in \cD$. We will show

\begin{theorem}
\label{T:bipara}
$[\bfT, b]$ is $L^2$-bounded if and only if $b \in bmo^d$.
\end{theorem}

Recall that the space $bmo^d$, or dyadic $little \: bmo$, is defined by
	$$\|b\|_{bmo^d} := \sup_{R\in \cD^2} \frac{1}{|R|} \int_R |b(x, y) - \La b\Ra_R|^2\,d(x,y).$$
In terms of Haar functions, this is the same as
	\begin{equation}
	\label{E:littlebmo-def}
	\|b\|_{bmo^d} = \sup_{R_0 = I_0\times J_0} \left(
		\frac{1}{|I_0||J_0|} \sum_{\substack{I\subseteq I_0\\ J \subseteq J_0}} |(b, h_I \otimes h_J)|^2
		+ \frac{1}{|I_0|} \sum_{I \subseteq I_0} |(b, h_I \otimes \tilde\one_{J_0})|^2
		+ \frac{1}{|J_0|} \sum_{J \subseteq J_0} |(b, \tilde\one_{I_0} \otimes h_J)|^2
	\right)
	\end{equation}
	
Above we had $D,  \pi, Z$ parts of the commutator. Now we will have the following  parts:
$$
DD, D\pi, \pi D, \pi\pi, \pi Z, Z\pi, DZ, ZD, ZZ,
$$
detailed below. As before, we consider $b$ fixed so the term $\bfT \pi\pi(b, f) - \pi\pi(b, \bfT f)$ will simply be denoted $\bfT\pi\pi - \pi\pi\bfT$, for example.
	
	\noindent\textbf{1.} $\pi\pi$ term:
	$$
	(\bfT\pi\pi - \pi\pi\bfT)f = \sum_{\substack{I even \\J even}} (b, h_I\otimes h_J) \La f\Ra_{I\times J} (T_1 h_I) \otimes (T_2 h_J) - \sum_{I, J} (b, h_I\otimes h_J) \La \bfT f\Ra_{I\times J} h_I \otimes h_J.
	$$

	\noindent\textbf{2.} $ZZ$ term:
	$$
	(\bfT ZZ - ZZ\bfT)f = \sum_{I, J} (b, h_I\otimes h_J) (f, h_I \otimes h_J) (T_1 \tilde\one_I)\otimes(T_2\tilde\one_J) - \sum_{\substack{I odd\\J odd}} (b, h_I\otimes h_J) (f, h_{\hat I}\otimes h_{\hat J}) 
	s(I, \hat I) s(J, \hat J) \tilde\one_I \otimes \tilde\one_J.
	$$

	\noindent\textbf{3.} $\pi Z$ term:
	$$
	(\bfT \pi Z - \pi Z\bfT)f = \sum_{\substack{I even \\ J}} (b, h_I\otimes h_J) (f, \tilde\one_I \otimes h_J)
	(T_1 h_I)\otimes(T_2\tilde\one_J) - \sum_{I, J} (b, h_I\otimes h_J) (\bfT f, \tilde\one_I\otimes h_J) h_I \otimes \tilde\one_J.
	$$

	\noindent\textbf{4.} $Z\pi$ term:
	$$
	(\bfT Z\pi - Z\pi\bfT)f = \sum_{\substack{I \\ J even}} (b, h_I\otimes h_J) (f, h_I\otimes\tilde\one_J) (T_1\tilde\one_I)\otimes(T_2 h_J) - \sum_{I, J} (b, h_I\otimes h_J) (\bfT f, h_I\otimes\tilde\one_J) \tilde\one_I\otimes h_J.
	$$

	\noindent\textbf{5.} $\pi D$ term:
	$$
	(\bfT\pi D - \pi D\bfT)f = \sum_{\substack{I even \\J even}} (b, h_I\otimes\tilde\one_J) (f, \tilde\one_I\otimes h_J) (T_1 h_I)\otimes(T_2 h_J) - \sum_{I,J} (b, h_I\otimes\tilde\one_J) (\bfT f, \tilde\one_I\otimes h_J) h_I\otimes h_J.
	$$
	
	\noindent\textbf{6.} $Z D$ term:
	$$
	(\bfT ZD - ZD\bfT)f = \sum_{\substack{I\\J even}} (b, h_I\otimes\tilde\one_J) (f, h_I\otimes h_J) (T_1\tilde\one_I)\otimes(T_2 h_J) - \sum_{\substack{I odd\\J odd}} (b, h_I\otimes\tilde\one_J) 
	(f, h_{\hat I}\otimes h_{\hat J}) s(I, \hat I) s(J, \hat J) \tilde\one_I\otimes h_J.
	$$	

	\noindent\textbf{7.} $D\pi$ term:
	$$
	(\bfT D\pi - D\pi\bfT)f = \sum_{\substack{I even\\J even}} (b, \tilde\one_I\otimes h_J) (f, h_I\otimes\tilde\one_J) (T_1 h_I)\otimes(T_2 h_J) - \sum_{I,J} (b, \tilde\one_I\otimes h_J) 
	(\bfT f, h_I\otimes\tilde\one_J) h_I\otimes h_J.
	$$

	\noindent\textbf{8.} $DZ$ term:
	$$
	(\bfT DZ - DZ\bfT)f = \sum_{\substack{I even\\J}} (b, \tilde\one_I\otimes h_J) (f, h_I \otimes h_J) 
	(T_1 h_I)\otimes(T_2 \tilde\one_J) - \sum_{\substack{I odd\\J odd}} (b, \tilde\one_I \otimes h_J)
	(f, h_{\hat I}\otimes h_{\hat J}) s(I, \hat I) s(J, \hat J) h_I\otimes\tilde\one_J.
	$$

	\noindent\textbf{9.} $DD$ term:
	$$
	(\bfT DD - DD\bfT)f = \sum_{\substack{I odd\\J odd}} \big( \La b\Ra_{\hat I \times \hat J}
		- \La b\Ra_{I\times J}\big) (f, h_{\hat I}\otimes h_{\hat J}) s(I, \hat I) s(J, \hat J) h_I\otimes h_J.
	$$

We also have a key lemma for this proof:
\begin{lemma}
\label{L:basebounds}
If $[\bfT, b]$ is bounded, then:
	\begin{equation}
	\label{E:basebd1}
	|(b, h_I \otimes h_J)|^2 \lesssim |I||J|,
	\end{equation}
	
	\begin{equation}
	\label{E:basebd2}
	|(b, \tilde\one_I \otimes h_J)|^2 \lesssim |J|,
	\end{equation}
	
	\begin{equation}
	\label{E:basebd3}
	|(b, h_I \otimes \tilde\one_J)|^2 \lesssim |I|,
	\end{equation}
for all $I, J \in \cD$.
\end{lemma}
Above and from here on out, the symbol $\lesssim$ means that the quantity is  bounded by some constant multiple of the $L^2$-norm of $[\bfT, b]$. We prove this lemma in section \ref{Ss:basebounds}.

Remark that this Lemma is the full bi-parameter equivalent of \eqref{dia}. Simply testing on $h_{I_0}\otimes h_{J_0}$, where both $I_0, J_0$ are even (the simple analog of the one-parameter case) will indeed give us that
	\begin{equation}
	\label{E:basedecoy}
	|\La b\Ra_{\hat I \times \hat J}
		- \La b\Ra_{I\times J}|\lesssim 1, \forall I, J \text{ both odd,}
	\end{equation}
which does bound the $DD$ term above -- as detailed in Section \ref{Ss:basebounds}. However, even having \eqref{E:basedecoy} for all $I,J$ is not enough for our purposes, because in two parameters
	$$
	\La b\Ra_{\hat I \times \hat J}
		- \La b\Ra_{I\times J} = (b, h_{\hat I}\otimes h_{\hat J}) h_{\hat I}(I) h_{\hat J}(J)
		+ (b, \tilde\one_{\hat I}\otimes h_{\hat J}) h_{\hat J}(J)
		+ (b, h_{\hat I} \otimes \tilde\one_{\hat J}) h_{\hat I}(I).
	$$
So Lemma \ref{L:basebounds} is much stronger.

%----------------proof lower bound 2par-----------------%
\subsection{Proof of Theorem \ref{T:bipara} - Lower Bound.}
\begin{proof}[\unskip\nopunct]
First of all, it is easy to see that Lemma \ref{L:basebounds} automatically bounds the operators $DD, ZZ, ZD$, and $DZ$ for all test functions. So these operators may be discarded. The remaining five operators, we now test on a Haar function $h_{I_0}\otimes h_{J_0}$. This yields ten operators, grouped by the various parity restrictions on $I_0$ and $J_0$:
	
	\vspace{0.1in}
	\noindent \textbf{Group 1:} Terms that require {\it both} $I_0$, $J_0$ be {\it even}:

\begin{fleqn}[\parindent]
\begin{eqnarray*}
\pi\pi2 &=& 
	- \sum_{\substack{I\subsetneq I_0^+ \\ J \subsetneq J_0^+}} (b, h_I\otimes h_J) h_{I_0^+}(I) h_{J_0^+}(J) h_I\otimes h_J 
	+ \sum_{\substack{I\subsetneq I_0^+ \\ J \subsetneq J_0^-}} (b, h_I\otimes h_J) h_{I_0^+}(I) h_{J_0^-}(J) h_I\otimes h_J \\
	&& 
	+ \sum_{\substack{I\subsetneq I_0^- \\ J \subsetneq J_0^+}} (b, h_I\otimes h_J) h_{I_0^-}(I) h_{J_0^+}(J) h_I\otimes h_J 
	- \sum_{\substack{I\subsetneq I_0^- \\ J \subsetneq J_0^-}} (b, h_I\otimes h_J) h_{I_0^-}(I) h_{J_0^-}(J) h_I\otimes h_J. 
\end{eqnarray*}
\end{fleqn}

\begin{fleqn}[\parindent]
\begin{eqnarray*}
\pi Z2 &=& 
	- \sum_{I \subsetneq I_0^+} (b, h_I\otimes h_{J_0^+}) h_{I_0^+}(I) h_I\otimes \tilde\one_{J_0^+}
	+ \sum_{I \subsetneq I_0^+} (b, h_I\otimes h_{J_0^-}) h_{I_0^+}(I) h_I\otimes \tilde\one_{J_0^-}\\
	&&+ \sum_{I \subsetneq I_0^-} (b, h_I\otimes h_{J_0^+}) h_{I_0^-}(I) h_I\otimes \tilde\one_{J_0^+}
	- \sum_{I \subsetneq I_0^-} (b, h_I\otimes h_{J_0^-}) h_{I_0^-}(I) h_I\otimes \tilde\one_{J_0^-}.
\end{eqnarray*}
\end{fleqn}

\begin{fleqn}[\parindent]
\begin{eqnarray*}
Z\pi2 &=& -\sum_{J\subsetneq J_0^+} (b, h_{I_0^+}\otimes h_J) h_{J_0^+}(J) \tilde\one_{I_0^+}\otimes h_J 
	+ \sum_{J\subsetneq J_0^+} (b, h_{I_0^-}\otimes h_J) h_{J_0^+}(J) \tilde\one_{I_0^-}\otimes h_J \\
	&& + \sum_{J\subsetneq J_0^-} (b, h_{I_0^+}\otimes h_J) h_{J_0^-}(J) \tilde\one_{I_0^+}\otimes h_J 
	-\sum_{J\subsetneq J_0^-} (b, h_{I_0^-}\otimes h_J) h_{J_0^-}(J) \tilde\one_{I_0^-}\otimes h_J. 
\end{eqnarray*}
\end{fleqn}

\begin{fleqn}[\parindent]
\begin{eqnarray*}
\pi D2 &=& - \sum_{I\subsetneq I_0^+} (b, h_I\otimes \tilde\one_{J_0^+}) h_{I_0^+}(I) h_I\otimes h_{J_0^+}
	+ \sum_{I\subsetneq I_0^+} (b, h_I\otimes \tilde\one_{J_0^-}) h_{I_0^+}(I) h_I\otimes h_{J_0^-}\\
	&& + \sum_{I\subsetneq I_0^-} (b, h_I\otimes \tilde\one_{J_0^+}) h_{I_0^-}(I) h_I\otimes h_{J_0^+}
	- \sum_{I\subsetneq I_0^-} (b, h_I\otimes \tilde\one_{J_0^-}) h_{I_0^-}(I) h_I\otimes h_{J_0^-}.
\end{eqnarray*}
\end{fleqn}

\begin{fleqn}[\parindent]
\begin{eqnarray*}
D\pi2 &=& - \sum_{J \subsetneq J_0^+} (b, \tilde\one_{I_0^+}\otimes h_J) h_{J_0^+}(J) h_{I_0^+}\otimes h_J
	+ \sum_{J \subsetneq J_0^+} (b, \tilde\one_{I_0^-}\otimes h_J) h_{J_0^+}(J) h_{I_0^-}\otimes h_J\\
	&&+ \sum_{J \subsetneq J_0^-} (b, \tilde\one_{I_0^+}\otimes h_J) h_{J_0^-}(J) h_{I_0^+}\otimes h_J
	- \sum_{J \subsetneq J_0^-} (b, \tilde\one_{I_0^-}\otimes h_J) h_{J_0^-}(J) h_{I_0^-}\otimes h_J.
\end{eqnarray*}
\end{fleqn}

	\vspace{0.1in}
	\noindent \textbf{Group 2:} Term that requires $I_0$ be {\it even}:
\begin{fleqn}[\parindent]
\begin{eqnarray*}
D\pi1 = \sum_{J even, \: J \subsetneq J_0} (b, \tilde\one_{I_0}\otimes h_J) h_{J_0}(J) (T_1 h_{I_0}) \otimes (T_2 h_J).
\end{eqnarray*}
\end{fleqn}

	\vspace{0.1in}
	\noindent \textbf{Group 3:} Term that requires $J_0$ be {\it even}:
\begin{fleqn}[\parindent]
\begin{eqnarray*}
\pi D1 = \sum_{I even, \: I \subsetneq I_0} (b, h_I\otimes \tilde\one_{J_0}) h_{I_0}(I) (T_1 h_I) \otimes (T_2 h_{J_0}).
\end{eqnarray*}
\end{fleqn}

	\vspace{0.1in}
	\noindent \textbf{Group 4:} Terms with no restrictions on the parity of $I_0$, $J_0$:
\begin{fleqn}[\parindent]
\begin{eqnarray*}
	\pi\pi1 = \sum_{\substack{I even,\:I\subsetneq I_0\\J even, \: J\subsetneq J_0}}
		(b, h_I\otimes h_J) h_{I_0}(I) h_{J_0}(J) (T_1h_I)\otimes(T_2h_J).
\end{eqnarray*}
\end{fleqn}

\begin{fleqn}[\parindent]
\begin{eqnarray*}
	\pi Z1 = \sum_{I even,\:I\subsetneq I_0}
		(b, h_I\otimes h_{J_0}) h_{I_0}(I)  (T_1h_I)\otimes(T_2\tilde\one_{J_0}).
\end{eqnarray*}
\end{fleqn}

\begin{fleqn}[\parindent]
\begin{eqnarray*}
	Z\pi1 = \sum_{J even, \: J\subsetneq J_0}
		(b, h_{I_0}\otimes h_J) h_{J_0}(J) (T_1\tilde\one_{I_0})\otimes(T_2h_J).
\end{eqnarray*}
\end{fleqn}

\vspace{0.1in}
\noindent \textbf{I.} {\it Take $I_0$, $J_0$ both odd}. In this case, we only have the three terms in Group 4 above, the terms with no restrictions on the parity of $I_0$, $J_0$. Note that all three of these terms are mutually orthogonal:
	\begin{eqnarray*}
	\pi\pi1 &:& (I\subsetneq I_0)\times(J \subsetneq J_0)\\
	\pi Z1 &:&  (I \subsetneq I_0)\times(J \supsetneq J_0)\\
	Z\pi1 &:& (I \supsetneq I_0) \times (J \subsetneq J_0),
	\end{eqnarray*}
so all three are {\it individually bounded} on the test function. Looking at the $L^2$-norm of $\pi\pi1$, for example, shows that
	$$
	\sum_{\substack{I even,\:I\subsetneq I_0\\J even, \: J\subsetneq J_0}}
		|(b, h_I\otimes h_J)|^2 \lesssim |I_0||J_0|, \text{ for all } I_0, J_0 \text{ both odd}.
	$$
But combining with Lemma \ref{L:basebounds}, we see that this is actually true for {\it all} $I_0$, $J_0$, meaning that $\pi\pi1$ is uniformly bounded on {\it all} test functions $h_{I_0}\otimes h_{J_0}$ and this term may be discarded. The same holds similarly for the remaining terms $\pi Z1$ and $Z\pi1$ in Group 4. With all Group 4 terms now being discarded, we move on.

\vspace{0.1in}
\noindent \textbf{II.} {\it Take $I_0$ even, $J_0$ odd.} In this case we only have the term $D\pi 1$ in Group 2. Computing the $L^2$-norm here gives
	$$
	\sum_{J even, J \subsetneq J_0} |(b, \tilde\one_{I_0}\otimes h_J)|^2 \lesssim |J_0|, \text{ for all odd } J_0.
	$$
Again, combined with Lemma \ref{L:basebounds}, this actually holds for all $J_0$, meaning that $D\pi1$ is uniformly bounded on all test functions $h_{I_0}\otimes h_{J_0}$, and may be discarded.

\vspace{0.1in}
\noindent \textbf{III.} {\it Take $I_0$ odd, $J_0$ even.} This case is symmetric to Case II, and yields that the term $\pi D1$ in Group 3 is uniformly bounded on all test functions $h_{I_0}\otimes h_{J_0}$, and may be discarded.

\vspace{0.1in}
\noindent \textbf{IV.} {\it Take $I_0$, $J_0$ both even.} The last case remaining, where we only need to look at the five terms in Group 1. But happily, these are also mutually orthogonal and therefore {\it individually bounded}!
	\begin{eqnarray*}
	\pi\pi2 : &(I\subsetneq I_0^\pm)&\times(J \subsetneq J_0^\pm)\\
	\pi Z2 :  &(I \subsetneq I_0^\pm)&\times(J \supsetneq J_0^\pm)\\
	Z\pi2 : &(I \supsetneq I_0^\pm) &\times (J \subsetneq J_0^\pm)\\
	\pi D2 : &(I \subsetneq I_0^\pm)&\times(J_0^\pm)\\
	Z\pi2 : &(I_0^\pm) &\times (J \subsetneq J_0^\pm)
	\end{eqnarray*}
Computing the $L^2$-norm of $\pi\pi2$ gives
	$$
	\sum_{I\subseteq I_0,\: J \subseteq J_0} |(b, h_I\otimes h_J)|^2 \lesssim |I_0||J_0|,
	$$
the $L^2$-norm of $\pi D2$ gives
	$$
	\sum_{I\subseteq I_0} |(b, h_I \otimes \tilde\one_{J_0})|^2 \lesssim |I_0|,
	$$
and finally the $L^2$-norm of $D\pi2$ gives
	$$
	\sum_{J\subseteq J_0} |(b, \tilde\one_{I_0}\otimes h_J)|^2 \lesssim |J_0|.
	$$	
These mean exactly that $b \in bmo^d$.
\end{proof}

%---------proof - lemma 2par-------%
\subsection{Proof of Lemma \ref{L:basebounds}.}
\label{Ss:basebounds}
\begin{proof}[\unskip\nopunct]
\noindent \textbf{I.} {\it Testing on $h_I\otimes h_J$.} If $[\bfT, b]$ is bounded, then $\|[\bfT, b]h_R\|\lesssim 1$ for all rectangles $R\in\cD^2$. In particular,
	$$
	\left| \left([\bfT, b]h_{I\times J}, \: h_{K\times L}\right)\right| \lesssim 1, \forall I, J, K, L \in\cD.
	$$
Further expanding this:
	\begin{equation}
	\label{E:base1}
	\left| \bigg(b h_I\otimes h_J, \: \underbrace{(T_1^* h_K)\otimes(T_2^*h_L)}_{0\text{ unless both }K, L\text{ odd}}
		\bigg)
	- \bigg(b \underbrace{(T_1 h_I)\otimes(T_2 h_J)}_{0 \text{ unless both } I,J \text{ even}}, h_K\otimes h_L\bigg)
	\right| \lesssim 1, \forall I, J, K, L \in \cD.
	\end{equation}

\vspace{0.1in}
\noindent \textbf{Ia.} Take {\it $I, J$ both even} and $K = I_{\pm}$, $L = J_{\pm}$ (children of $I, J$) -- so both $K, L$ are odd. Then both terms of \eqref{E:base1} are present, and it becomes exactly \eqref{E:basedecoy}.

\vspace{0.1in}
\noindent \textbf{Ib.} Take {\it $I, J$ both odd} (so second term in \eqref{E:base1} is $0$), and $K = I^{(2)}$, 
$L = J^{(2)}$. Both $K, L$ are then odd, and $T_1^* h_K = \pm h_{I^{(3)}}$, $T_2^* h_L = \pm h_{J^{(3)}}$. 
Equation \eqref{E:base1} then becomes
	$$
	\left|\bigg( b h_I\otimes h_J, \: h_{I^{(3)}}\otimes h_{J^{(3)}} \bigg)\right| \lesssim 1.
	$$
But since $h_{I^{(3)}}$, $h_{J^{(3)}}$ are constant on $I$, $J$, respectively:
	$$h_{I^{(3)}}(I) = \pm \frac{1}{\sqrt{|I^{(3)}|}}; \:\:\:\:
	h_{J^{(3)}}(J) = \pm \frac{1}{\sqrt{|J^{(3)}|}},$$
we now have \eqref{E:basebd1} for the case when $I, J$ are both odd:
	$$
	|(b, h_I \otimes h_J)|^2 \lesssim |I||J|, \forall I, J \text{ both odd}.
	$$
	
\vspace{0.1in}
\noindent \textbf{Ic.} Take {\it $I$ even, $J$ odd} (so second term in \eqref{E:base1} is $0$), and
$K = I_{\pm}$, $L = J^{(2)}$ (both odd). Then $T_1^* h_K = \pm h_I$, and $T_2^* h_L = \pm h_{J^{(3)}}$.
Equation \eqref{E:base1} becomes
	$$
	\left|\bigg( b h_I\otimes h_J, \: h_I\otimes h_{J^{(3)}} \bigg)\right| \lesssim 1,
	$$
which yields \eqref{E:basebd2} for the case when $I$ is even and $J$ is odd:
	$$
	|(b, \tilde\one_I \otimes h_J)|^2 \lesssim |J|, \forall I \text{ even, } J \text{ odd.}
	$$
	
\vspace{0.1in}
\noindent \textbf{Id.} Take {\it $I$ odd, $J$ even} (so second term in \eqref{E:base1} is $0$), and
$K = I^{(2)}$, $L = J_{\pm}$ (both odd). Symmetrically with case Ic, this will give \eqref{E:basebd3}
for the case of $I$ odd and $J$ even:
	$$
	|(b, h_I \otimes \tilde\one_J)|^2 \lesssim |I|, \forall I \text { odd, } J \text{ even.}
	$$

\vspace{0.1in}
\noindent \textbf{Ie.}	Take {\it $I$, $J$ both even}, and $K=I$ (even), $L = J_{\pm}$ (odd). Then the first term in \eqref{E:base1} is $0$, and we have
	$$
	\left| \bigg( b(h_{I_+} - h_{I_-})\otimes(h_{J_+} - h_{J_-}), h_I \otimes h_{J_\pm}
	\bigg)\right| \lesssim 1.
	$$
This gives us \eqref{E:basebd3} for the case when $I, J$ are both odd. 

\vspace{0.1in}
\noindent \textbf{If.} Symmetrically, take {\it $I$, $J$ both even}, $K=I_{\pm}$ (odd), $L = J$ (even), and
we have \eqref{E:basebd2} for the case when $I, J$ are both odd.

\vspace{0.2in}
\noindent \textbf{II.} {\it Testing on $h_I\otimes T_2^*h_J$.} Running through the same process with the test function $h_I\otimes T_2^*h_J$ instead of $h_I\otimes h_J$, we get
	\begin{equation}
	\label{E:base2}
	\left| \bigg(\underbrace{b h_I\otimes T_2^*h_J}_{0 \text{ unless } J \text{ odd}}, \: \underbrace{(T_1^* h_K)\otimes(T_2^*h_L)}_{0\text{ unless both }K, L\text{ odd}}
		\bigg)
	- \bigg(b \underbrace{(T_1 h_I)\otimes(T_2 T_2^* h_J)}_{0 \text{ unless } I\text{ even, } J \text{ odd}}, h_K\otimes h_L\bigg)
	\right| \lesssim 1, \forall I, J, K, L \in \cD.
	\end{equation}
Note that, in one parameter,
	$$
	T T^* h_I = h_I - h_{s(I)}, \forall I \text{ odd; } 0 \text{ otherwise.}
	$$
where $s(I)$ denotes the dyadic sibling of $I$.

\vspace{0.1in}
\noindent \textbf{IIa.} Take $I$ even, $J$ odd, and $K = \hat I$ (odd), $L = J$ (odd). Then both terms of \eqref{E:base2} are present, and give 
	$$
	\left|\frac{1}{|I^{(2)}|} (b, h_I \otimes \tilde\one_{\hat J}) 
	- s(I, \hat I) \frac{1}{|\hat I|} \bigg( b, \: (h_{I_+} - h_{I_-})\otimes\tilde\one_J \bigg)\right| \lesssim 1.
	$$
Since $I_{\pm}$ and $J$ are odd, we already know from case {\bf Ie} that
	$$
	\left|s(I, \hat I) \frac{1}{|\hat I|} 
	\bigg( b, \: (h_{I_+} - h_{I_-})\otimes\tilde\one_J \bigg)\right| \lesssim 1,
	$$
so we now have \eqref{E:basebd3} for $I, J$ both even.

\vspace{0.1in}
\noindent \textbf{IIb.} Take $I$ even, $J$ odd, and $K = I_\pm$ (odd), $L = J_\pm$ (even). Then the first term in \eqref{E:base2} is $0$, and the second term is
	$$
	\left|\bigg(
		b(h_{I_+} - h_{I_-})\otimes(h_J - h_{s(J)}), \: h_{I_\pm} \otimes h_{J_\pm}
	\bigg)\right| = \frac{1}{\sqrt{|J|}} \left|
	\bigg(b, \tilde\one_{I_\pm} \otimes h_{J_\pm}\bigg) \right| \lesssim 1,
	$$
so now we have \eqref{E:basebd2} for the case where $I$ is odd and $J$ is even.

\vspace{0.1in}
\noindent \textbf{IIc.} Take $I$ even, $J$ odd, and $K = I$ (even), $L = J_\pm$ (even). Again the first term in \eqref{E:base2} is $0$, and we have
	$$
	\left|\bigg(
		b(h_{I_+} - h_{I_-})\otimes(h_J - h_{s(J)}), \: h_I \otimes h_{J_\pm}
	\bigg)\right| = \frac{1}{\sqrt{|I|}\sqrt{|J|}} \left|
	\bigg(b, h_{I_\pm} \otimes h_{J_\pm}\bigg) \right| \lesssim 1,
	$$
which gives us \eqref{E:basebd1} for $I$ odd and $J$ even.

\vspace{0.2in}
\noindent \textbf{III.} {\it Testing on $T_1^*h_I\otimes h_J$.} This case is symmetrical to case II, and it gives us \eqref{E:basebd2} for $I, J$ both even, \eqref{E:basebd3} for $I$ even and $J$ odd, and \eqref{E:basebd1} for $I$ even and $J$ odd. 

At this point \eqref{E:basebd2} and \eqref{E:basebd3} are fully proved, and all that is left is \eqref{E:basebd1} for $I, J$ both even. This final situation is dealt with below.

\vspace{0.2in}
\noindent \textbf{IV.} {\it Testing on $T_1^*h_I\otimes T_2^*h_J$.} In this case we have
$$
	\left| \bigg(\underbrace{b (T_1^*h_I)\otimes (T_2^*h_J)}_{0 \text{ unless } I, J \text{ both odd}}, \: \underbrace{(T_1^* h_K)\otimes(T_2^*h_L)}_{0\text{ unless both }K, L\text{ odd}}
		\bigg)
	- \bigg(b \underbrace{(T_1 T_1^* h_I)\otimes(T_2 T_2^* h_J)}_{0 \text{ unless } I, J\text{ both odd}}, h_K\otimes h_L\bigg)
	\right| \lesssim 1, \forall I, J, K, L \in \cD.
$$
In the equation above, take $I, J$ both odd, and $K = I_\pm$, $L = J_\pm$ (both even), so we only have the second term above:
	$$
	\left|\bigg(
		b(h_I - h_{s(I)})\otimes(h_J - h_{s(J)}), \: h_{I_\pm} \otimes h_{J_\pm}
	\bigg)\right| = \frac{1}{\sqrt{|I|}\sqrt{|J|}} \left|
	\bigg(b, h_{I_\pm} \otimes h_{J_\pm}\bigg) \right| \lesssim 1,
	$$
which gives us the last piece: \eqref{E:basebd1} for $I, J$ both even.
\end{proof}

%------------proofs upper bound 2par-------------%
\subsection{Proof of Theorem \ref{T:bipara} - Upper Bound.}
\begin{proof}[\unskip\nopunct]
Again, the proof that $[\bfT, b]$ is bounded if $b\in bmo^d$ is completely standard, but we include it here for completeness.
We will show that each term in the paraproduct decomposition of $[\bfT, b]$ is bounded. Some of these will pair with {\it product} BMO and others with {\it little} $bmo$. The easiest one is again the $DD$ term, since its $L^2$-norm can simply be computed:
	$$
	\|(\bfT DD - DD\bfT)\|_2^2 \leq \sum_{I, J} \left|\La b\Ra_{\hat I\times\hat J} - \La b\Ra_{I\times J}\right|^2
		|(f, h_{\hat I}\otimes h_{\hat J})|^2 \lesssim \|b\|_{bmo^d} \|f\|_2^2.
	$$
Of the remaining eight terms, four will ``pair'' with {\it product} $BMO^d$, and four will ``pair'' with $bmo^d$. 

\subsubsection{Product $BMO^d$ terms: $\pi\pi$, $ZZ$, $\pi Z$, $Z\pi$.} For these terms, we use the fact that if $b \in bmo^d$ then $b \in BMO^d$, and the bi-parameter $H_1^d-BMO^d$ duality:
	$
	|(b, \Phi)| \lesssim \|b\|_{BMO^d} \|S_d \Phi\|_1,
	$
where $S_d$ now denotes the bi-parameter dyadic square function:
	$
	S_d^2 f = \sum_{I, J} |(f, h_I\otimes h_J)|^2\tilde\one_I\otimes\tilde\one_J.
	$

We split the $(\bfT\pi\pi - \pi\pi\bfT)f$ term into $\Pi_1 f - \Pi_2f$, where
	$$
	\Pi_1f = \sum_{I,J} (b, h_I\otimes h_J) \La f\Ra_{I\times J} \bfT(h_I \otimes h_J); \:\:\:\:
	\Pi_2f = \sum_{I, J} (b, h_I\otimes h_J)\La \bfT f\Ra_{I\times J} h_I\otimes h_J.
	$$
For the first term:
	$
	(\Pi_1f, g) = (b, \Phi_1)
	$, where
	$
	\Phi_1 = \sum_{I, J} \La f\Ra_{I\times J} (\bfT^*g, h_I\otimes h_J) h_I\otimes h_J,
	$
and $f, g\in L^2(\bR^2)$. Then
	$$
	S_d^2\Phi_1 = \sum_{I,J} \La f\Ra_{I\times J}^2 |(\bfT^*g, h_I\otimes h_J)|^2 \tilde\one_{I\times J}
	\leq [M_d^2f] \: [S_d^2(T^*g)],
	$$
where $M_d$ now denotes the bi-parameter dyadic maximal function. Then
	$$
	|(\Pi_1f, g)| = |(b, \Phi_1)| \lesssim \|b\|_{BMO^d} \|S_d\Phi_1\|_1 \lesssim \|b\|_{BMO^d}\|M_df\|_2\|S_d \bfT^*g\|_2
	\lesssim  \|b\|_{BMO^d}\|f\|_2\|g\|_2.
	$$
The other term $\Pi_2$, as well as the terms in $(\bfT ZZ - ZZ\bfT)$ follow similarly.

For the $\pi Z$ and $Z\pi$ terms, we appeal to the mixed square functions introduced in \cite{HPW}:
	$$
	[SM]^2f(x, y) := \sum_I M^2_{d_2}f_I(y) \tilde\one_I(x), \:\:\:\:
	[MS]^2f(x, y) := \sum_J M^2_{d_1}f_J(x) \tilde\one_J(y),
	$$
where $M_{d_i}$ denote the dyadic maximal function in parameter $i$ only, and for every $I, J$:
	$$
	f_I(y) := \int_\bR f(x, y) h_I(x)\,dx; \:\:\:\:
	f_J(x) := \int_\bR f(x, y) h_J(y)\,dy.
	$$
More general forms of these mixed functions were proved to be bounded in the weighted setting in \cite{HPW}. 

Now, looking at our term $(\bfT\pi Z - \pi Z\bfT)f = \Pi Z_1f - \Pi Z_2f$, where
	$$
	\Pi Z_1 f = \sum_{\substack{I even \\ J}} (b, h_I\otimes h_J) (f, \tilde\one_I \otimes h_J) \bfT(h_I\otimes\tilde\one_J); 
	\:\:\:\:
	\Pi Z_2 f = \sum_{I, J} (b, h_I\otimes h_J) (\bfT f, \tilde\one_I\otimes h_J) h_I \otimes \tilde\one_J.
	$$
For the first term, $(\Pi Z_1 f, g) = (b, \Phi_1)$, where
	$
	\Phi_1 = \sum_{I, J} \La f_J\Ra_I \La (\bfT^*g)_I \Ra_J	h_I\otimes h_J.
	$
Then
	\begin{eqnarray*}
	S_d^2\Phi_1 &=& \sum_{I, J} |\La f_J\Ra|^2 |\La (\bfT^*g)_I \Ra_J|^2 \tilde\one_I\otimes\tilde\one_J \\
	&\leq& \bigg( \sum_I M^2_{d_2}(\bfT^*g)_I(y) \tilde\one_I(x)\bigg)
	\bigg( \sum_J M^2_{d_1}f_J(x) \tilde\one_J(y)\bigg) = [SM]^2(\bfT^*g) [MS]^2(f).
	\end{eqnarray*}
So 
	$$
	|(\Pi Z_1 f, g)| \lesssim \|b\|_{BMO^d} \|S_d\Phi_1\|_1 \lesssim \|b\|_{BMO^d} \|[SM]^2(\bfT^*g)\|_2 \|[MS]^2(f)\|_2
	\lesssim \|b\|_{BMO^d}\|f\|_2\|g\|_2.
	$$
The other term $\Pi Z_2$, as well as the terms in $(Z\pi-\pi Z)$ follow similarly.

\vspace{0.1in}
\subsubsection{Little $bmo^d$ terms: $\pi D$, $D\pi$, $ZD$, $DZ$.} For these terms we appeal to the dyadic square functions in one-parameter
	$$
	S_{d_1}^2f(x, y) = \sum_I f_I^2(y)\tilde\one_I(x); \:\:\:\:
	S_{d_2}^2f(x, y) = \sum_J f_J^2(x)\tilde\one_J(y),
	$$
and the duality with $bmo^d$:
	$$
	|(b, \Phi)| \lesssim \|b\|_{bmo^d} \|S_{d_i}\Phi\|_1.
	$$
These were all treated in the weighted situation in \cite{HPW}. For example for the $\pi D$ term we have
	$(\bfT\pi D - \pi D\bfT)f = \pi D_1f - \pi D_2 f$ where
	$$
	\pi D_1f = \sum_{I,J} (b, h_I\otimes\tilde\one_J) (f, \tilde\one_I\otimes h_J) \bfT(h_I\otimes h_J);
	\:\:\:\:
	\pi D_1f = \sum_{I,J} (b, h_I\otimes\tilde\one_J) (\bfT f, \tilde\one_I\otimes h_J) h_I\otimes h_J.
	$$
For the first term, $(\pi D_1f, g) = (b, \Phi_1)$ where
	$
	\Phi_1 = \sum_{I,J} (f, \tilde\one_I\otimes h_J) (\bfT^*g, h_I\otimes h_J) h_I*\tilde\one_J.
	$
Then
	\begin{eqnarray*}
	S_{d_1}^2\Phi_1 &=& \sum_I \bigg(\sum_J \La f_J\Ra_I (\bfT^*g, h_I\otimes h_J)\tilde\one_J\bigg)^2\tilde\one_I\\
	&\leq& \sum_I \bigg(\sum_J \La |f_J|\Ra_I^2 \tilde\one_J\bigg) 
		\bigg(\sum_J (\bfT^*g, h_I\otimes h_J)^2\tilde\one_J\bigg)\tilde\one_I\\
	&\leq& \sum_J M^2_{d_1}f_J(x)\tilde\one_J(y) \cdot \sum_{I,J} |(\bfT^*g, h_I\otimes h_J)|^2 \tilde\one_I\otimes\tilde\one_J
	= [MS]^2f \cdot S_d^2(\bfT^*g).
	\end{eqnarray*}
So 
	$$
	|(\pi D_1f, g)| \lesssim \|b\|_{bmo^d} \|S_{d_1}\Phi_1\|_1 \leq \|b\|_{bmo^d} \|[MS]f \cdot S_d(\bfT^*g)\|_1 
	\lesssim \|b\|_{bmo^d} \|f\|_2\|g\|_2.
	$$
The other term $\pi Z_2$, as well as the terms of $(\bfT D\pi-D\pi\bfT)$, $(\bfT ZD - ZD\bfT)$ and $(\bfT DZ - DZ\bfT)$ follow similarly.
\end{proof}

%%%%%%%%%%%%%%%%%%%%%%%%%%%%%%%%%%%%%%%%

\end{document}